\documentclass[12pt]{article}
\usepackage{amssymb}
\usepackage{amsmath, amscd}
\usepackage{diagrams}

\usepackage{mathptmx}

\newtheorem{prop}{Proposition}

\newtheorem{thm}{Theorem}

\newcommand{\h}{\hom }
\newcommand{\LL}{\mathbb L}

\begin{document}
\begin{center}
    \large{\bf Geometric algebra of projective lines}\\
    \normalsize
    \medskip
    Anders Kock\\
    University of Aarhus\\
    kock{@}imf.au.dk

    \end{center}
    
    \small
    
    \noindent Abstract. The projective line over a field carries  
      structure of a groupoid with a certain correspondence between 
      objects and arrows. We discuss to what extent the field can be 
      reconstructed from the groupoid.
      \normalsize
      
      \bigskip
    
    We consider a transitive groupoid $\LL$ where for any two 
    different objects $A$ and $B$, there is given a bijective 
    correspondence between the set $\h (A,B)$ and the set of objects 
    $C$ with $C\neq A$ and $C\neq B$. We say that $C$ is {\em label} 
    for the arrow $A\to B$. Note that  $\LL$ has at least 
    three objects (unless it is empty or a group; we exclude these 
    cases).
    
    Transitive groupoids with this structure, we call {\em projective 
    line candidates}.
     If $K$ is a field, the projective line $P(K^{2})$ over $K$ 
     gives rise to such structure (cf.\ \cite{CAPS}, \cite{APL}): the objects are the points of 
     $P(K^{2})$, i.e.\ the 1-dimensional linear subspaces of $K^{2}$; 
     the arrows are the linear isomorphisms between such. Such linear 
     isomorphism $A\to B$, for $A\neq B$, is projection in the 
     direction of a specific 
     direction $C$ with $C\neq A, B$, and this determines a bijective 
     correspondence of the kind postulated by the ``projective line candidate'' 
     notion.
     
     The  result of the present note is that {\em if a projective line candidate enjoys 
     certain properties of geometric character, then it is  
     isomorphic to a  $P(K^{2})$}, with $K$ an essentially unique 
     field.  (The notion of {\em isomorphism} of projective 
     line candidates is evident: an isomorphism of groupoids, 
     compatible with the assumed labelling of arrows  by objects.)
     
    \medskip
    
    \noindent{\bf Remark.} It is possible for the present purpose to replace the basic 
    structure of groupoid with a weaker notion of ``near-groupoid'', 
    which is like a groupoid, except that no endo-arrows $A\to A$ are 
    assumed. This reduction has here the advantage  that {\em all} 
    non-endo arrows in  $\LL$ ``are geometric'', i.e.\ can be {\em drawn} as 
    actual projections (or ``perspectivities'') , as in \cite{CAPS}, 
    \cite{DL}, \cite{APL}. An example of such a geometric picture is 
    given in the Appendix.

    \medskip
    
    If $C$ is 
     label for an arrow $A\to B$,  we write 
     $$C:A\to B\quad \mbox{ or }\quad \begin{diagram}A&\rTo 
     ^{C}&B\end{diagram}$$ and 
     similar  standard diagrammatic 
     notation. 
     We compose from left to right. 
     
     The following is an attempt to {\em 
     construct} the field $K$ from the groupoid, essentially by 
     providing the set $\h (A,A) \cup \{0\}$ (for some arbitrarily 
     chosen\footnote{or more canonically, take $K$ as a quotient of the disjoint 
     union of all the $\h (A,A) \cup \{0\}$ as $A$ ranges over all 
     objects}  object $A$) with the structure of 
     field, retaining as multiplication the composition operation on 
     $\h (A,A)$, as given 
     by the groupoid structure. This is also the reason we shall talk 
     about arrows $A\to A$ as {\em scalars} (another reasonable 
     terminology is: {\em pure quantities}, cf.\ \cite{MSPQ}).
     
     \medskip
     
     We assume the following basic equations (other assumptions will 
     be called for later):
     \begin{align}\begin{CD}(A@>C>>B@>C>>A ) \end{CD}&= 1_{A}
	\label{one}\\
    \begin{CD}(A@>C>>B@>C>>D )\end{CD} &= \begin{CD}(A@>C>>D)\end{CD}
	\label{two}\end{align}
	
	\medskip
    We assume that the vertex groups $\h (A,A)$ of $\LL$ are 
    commutative. (This can be stated in a way which does not involve 
    any endo-arrow, i.e.\ it can be stated as a property of 
    near-groupoids, in the sense of the Remark above; namely, for any 
    three parallel arrows $f_{i}:A\to B$ (with $A\neq B$),
    $$f_{1}\cdot f_{2}^{-1}\cdot f_{3}= f_{3}\cdot f_{2}^{-1}\cdot 
    f_{1}:A\to B.$$ If one draws these two three-fold composites in 
    the projective line $P(K^{2})$, the geometric figure that arises 
    is the Pappus configuration. So for a projective line candidate, 
    commutativity = validity of Pappus' Axiom.)

    So $\h (A,A)$ 
    is {\em canonically} isomorphic to $\h (B,B)$, by conjugation by 
    some, hence any, $A\to B$ (arrows $A\to B$ exist, since we 
    assume $\LL$ transitive).

    If $\mu \in \h (A,A)$, we say that $\mu$ is a {\em scalar at $A$}; if 
    $\mu ' \in \h (B,B)$ corresponds to it under the  
    conjugation correspondence, we write $\mu  \equiv \mu '$.
    
    Consider four objects $A,B,C,D$, with $A,B,C$ mutually distinct, 
    and $A,B,D$ mutually distinct. If $\begin{diagram} A&\rTo ^{C}& B
    \end{diagram}$ and $\begin{diagram} A&\rTo ^{D}& B
    \end{diagram}$, we write $(A,B;C,D)$ or $\left[ 
    \begin{array}{cc}A&B\\
	C&D
    \end{array}\right]$
    for the scalar at $A$ given 
    as the composite
    $$\begin{diagram} A&\rTo ^{C}& B&\rTo ^{D}& A
    \end{diagram};$$ this is the classical {\em cross ratio} or {\em 
    bi-rapport} of $A,B,C,D$, see \cite{CAPS}, \cite{DL}, \cite{APL}. 
    The cross ratio is thus,  in display, the composite
    $$
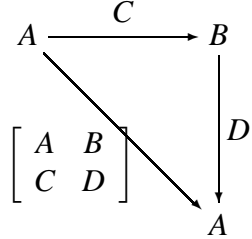
\begin{diagram}[nohug]A&\rTo ^{C}& B\\
    &\rdTo _{\left[ \begin{array}{cc}A&B\\C&D
    \end{array}\right] }&\dTo _{D}\\
    && A
    \end{diagram}$$
    Note that $(A,B;C,C) =1$, by (\ref{one}).
    The columns in the matrix displayed may be interchanged, modulo $\equiv$; 
    for, consider the diagram
    $$
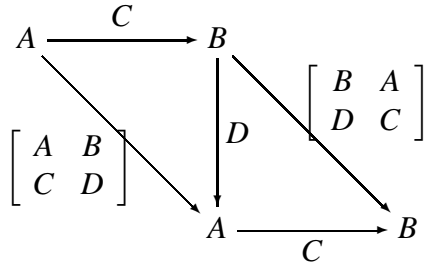
\begin{diagram}[nohug]A&\rTo ^{C}&B&&\\
    &\rdTo _{\left[ \begin{array}{cc}A&B\\C&D
    \end{array}\right] }&\dTo _{D}&
    \rdTo ^{\left[ \begin{array}{cc}B&A\\D&C
    \end{array}\right] }&\\
    &&A&\rTo_{C}&B
    \end{diagram}$$
    The two triangles commute by definition; so the commutativity of 
    the total quadrangle says that $C:A\to B$ conjugates the scalar
    $(A,B;C,D)$ at $A$ to the scalar $(B,A;D,C)$ at $B$,
    or
    $$\left[ \begin{array}{cc}A&B\\C&D
    \end{array}\right]\equiv \left[ \begin{array}{cc}B&A\\D&C
    \end{array}\right].$$
    The interchange of rows is less trivial, since it involves 
    change of names of objects into labels for arrows; we 
    impose as an axiom, the ``hexagon axiom'', which says that, for 
    $A,B,C,D$ mutually distinct, the (outer) 
    diagram
    
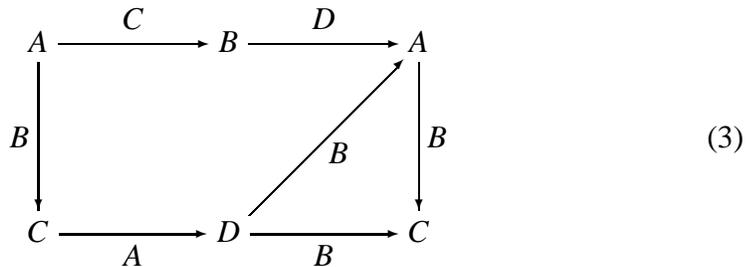
\begin{equation}\label{hex1}\begin{diagram}[nohug]A&\rTo ^{C}&B&\rTo ^{D}&A\\
    \dTo ^{B}&&&\ruTo _{B}&\dTo _{B}\\
    C&\rTo _{A}&D&\rTo _{B}&C
    \end{diagram}\end{equation}
    commutes, thus $B:A\to C$ conjugates the scalar $(A,B;C,D)$ at $A$ 
    to the scalar $(C,D;A,B)$ at $C$, i.e.
    $$\left[ \begin{array}{cc}A&B\\C&D
    \end{array}\right]\equiv \left[ \begin{array}{cc}C&D\\A&B
    \end{array}\right].$$
    Note that the two vertical arrows both are $B:A\to C$, but they may be 
    replaced jointly by any other $X:A\to C$ -- the conjugation relation is 
    the same, by commutativity of the groups of scalars. But note 
    that the advantage of putting $B$ is that then the hexagon can be 
    reduced to a pentagon (the ``inner'' pentagon in the diagram), because  two arrows with label $B$ in the right hand of 
    the diagram may be replaced by one single $B:D\to C$, by Axiom 
    (\ref{two}). Similarly if we put the label $D$ on both vertical 
    arrows. 
    
    It follows that cross ratios are invariant under the action of 
    the four-group, which is the rationale for the ``matrix'' 
    notation employed. Therefore also, any cross ratio in which the 
    letter $A$ occurs, may be replaced (up to $\equiv$) by one in 
    which the letter $A$ occurs in the the upper left hand corner, 
    without changing the cross ratio (mod $\equiv$).
    
    Note that $(A,B;C,D)=(A,B;C,D')$ implies $D=D'$.
    
    \medskip
    
    We now consider the effects of permuting the four vertices by a 
    permutation which is not one of the four-group permutations. 
    Since one of the letters, say $A$, may always be brought to the 
    upper left corner, it suffices to consider the permutations of the 
    remaining three entries $B,C,D$ (assumed mutually distinct, and distinct 
    from $A$).
    For cross ratios in 
    projective lines over a field, ghis gives a classical list  (cf.\ e.g. \cite{Struik} I-4);  
    it is reproduced here:
    let $\mu$ denote $(A,B;C,D)$ (this is recorded as the first equation in the 
    list). Then
    \begin{align}\left[ \begin{array}{cc}A&B\\C&D
    \end{array}\right]&= \mu \label{BCD}\\
    \left[ \begin{array}{cc}A&B\\D&C
    \end{array}\right]&=\mu ^{-1}\label{BDC}\\
    \left[ \begin{array}{cc}A&C\\B&D
    \end{array}\right]&=1-\mu\label{CBD}\\
    \left[ \begin{array}{cc}A&C\\D&B 
    \end{array}\right]&=(1-\mu )^{-1}\label{CDB}\\
    \left[ \begin{array}{cc}A&D\\B&C
    \end{array}\right] &= 1- \mu ^{-1}\label{DBC}\\
    \left[ \begin{array}{cc}A&D\\C&B
    \end{array}\right]&=(1-\mu ^{-1})^{-1}\label{DCB}
    \end{align}
    
    Equation (\ref{BDC}) makes sense and is easy to prove in our 
    context, using (\ref{one}). But the rest make 
    no sense as they stand, because we have not assumed any further 
    algebraic structure on the vertex group $\h (A,A)$ to justify the minus signs. 
    The crucial point is to give meaning to the right hand side of  (\ref{CBD}); the rest follow by combining 
    (\ref{BDC}) and (\ref{CBD}). What is meant by $1-\mu$?
    
    What is true for projective lines over a field is the following 
    property of cross ratios:
    \begin{equation}\mbox{{\em if $(A,B;C,D) = (A',B';C',D')$, then also $(A,C;B,D)= 
    (A',C';B',D')$}}\label{AS}\end{equation}
    To the extent this holds in $\LL$, we may define 
    an involution $\Phi$ on $\h (A,A)$: to define $\Phi (\mu )$ for a 
    scalar $\mu :A\to A$, we 
    {\em choose} $B,C,D$ so that 
    $$\mu =\left[ \begin{array}{cc}A&B\\C&D
    \end{array}\right],$$
    and then we put
    $$\Phi (\mu ) = \left[ \begin{array}{cc}A&C\\B&D
    \end{array}\right].$$
    For, by the  property assumed, the result $\Phi (\mu )$ does not 
    depend on the way $B,C,D$ were chosen. Also, we can then prove 
    (using variation of $A$) that the $\Phi$s thus defined on each 
    $\h (A,A)$ is invariant under the 
    (conjugation-) identification of $\h (A,A)$ with $\h (A',A')$.
    
    Under these circumstances, there is no harm in denoting $\Phi 
    (\mu )$ by $1-\mu $, and this we shall do.
    
    \medskip
    
    (So we assume the  property (\ref{AS}) as an axiom, but it is 
    unfortunately not purely equational, which we would prefer. I am 
    still looking for an equational formulation.)
    
    \medskip

    There is another {\em unary} ``minus'' operation possible, uniformly 
    on all the $\h (A,A)$s. We put
     $-\mu := (-1)\cdot \mu $, where $(-1):A\to A$ is the scalar at $A$ defined 
     as follows. We {\em choose} $B$ and $C$ (distinct, and distinct from 
     $A$) and let $(-1)_{A}$ be the scalar at $A$ defined as the 
     composite
     $$\begin{diagram}A&\rTo^{C} &B & 
     \rTo^{A}&C&\rTo^{B}&A\end{diagram};$$
      this particular composite is in the 
     coordinate situation (or in a projective line embedded in a 
     projective plane) (multiplication by) the scalar $-1$, see 
     \cite{CAPS}. It cannot be reduced to a 
      cross ratio, and it is a special case of 
     composites, considered in \cite{DL} under the name ``tri-rapport'' 
     (where cross ratio = ``bi-rapport''), see below.
     
     The independence   of $-1$ on the choice of $A,B,C$ again seems to be 
     something we need to impose as an axiom;  here, it will follow 
     from a purely equational one, namely commutativity, for all 
     $A,B,C,B',C'$, of the hexagon
     \begin{equation}\begin{diagram}A&\rTo^{C}&B&\rTo^{A}&C\\
     \dTo^{C'}&&&&\dTo_{B}\\
     B'&\rTo_{A}&C'&\rTo_{B'}&A
     \end{diagram}\label{hex2}\end{equation}
     From this follows that for each $A$, we have a well defined 
     scalar $(-1)_{A}$ (independent of the choice of $B$ and $C$).
     Then it easily follows that $(-1)_{A}\equiv (-1)_{B}$. For, we 
     may pick $C$ so that $(-1)_{A}$ is represented, as above, by 
     $A,B,C$, and $(-1)_{B}$ is similarly represented by $B,C,A$. But 
     then $C:A\to B$ conjugates the chosen expression for $(-1)_{A}$ 
     to the chosen one for $(-1)_{B}$. In this sense, $-1$ is a 
     ``uniform'' scalar.
     
     \medskip
     
     Let us in (\ref{hex2}) take $C':=B$ and $B':=C$. 
     Then we get the equality
     $$\begin{diagram}A&\rTo^{C}&B&\rTo^{A}&C&\rTo^{B}&A
     \end{diagram}$$
     $$=\begin{diagram}A&\rTo^{B}&C&\rTo^{A}&B&\rTo^{C}&A
     \end{diagram} $$
     and therefore that
     $$\begin{CD}A@>C>> B@>A>>C@>B>>A@>C>>B@>A>>C@>B>>A\end{CD}$$
     equals
     $$\begin{CD}A@>B>>C@>A>>B@>C>>A@>C>>B@>A>>C@>B>>A\end{CD}$$
     and this reduces to $1_{A}$ by three applications of 
     (\ref{one}). Thus $(-1).(-1)=1$, partly justifying the notation.
     Let us record this:
     \begin{prop}\label{xx}The scalar $-1$ has the property that $(-1)\cdot 
	 (-1)=1$.
	 \end{prop}
     
     \medskip
     
     We shall consider the notion of {\em tri-rapports}, in analogy 
     with cross ratios, which are also classically called 
     bi-rapports. 
     
     For $A,B,C$ mutually distinct, and $D\neq 
     A,B$, $E\neq B,C$, $F\neq A,C$, we have a scalar at $A$ given the 
     composite
     $$\begin{diagram}A&\rTo^{D}&B&\rTo^{E}&C&\rTo^{F}&A
     \end{diagram};$$
 we denote it $(A,B,C;D,E,F)$ or $\left[\begin{array}{ccc} 
A&B&C\\D&E&F \end{array}\right].$ Note that $-1$ is 
such a tri-rapport, $(-1)_{A} =(A,B,C;C,A,B)$, cf.\  
\cite{CAPS}.

     \medskip

	We  have
     $$\left[ \begin{array}{ccc}A&B&C\\
     D&E&F \end{array}\right] \equiv \left[ \begin{array}{ccc}B&C&A\\
     E&F&D \end{array}\right]$$
     (cyclic permutation of columns); this is clear: $D:A\to B$ will 
     conjugate the composite defining the left hand side to the one 
     defining the right hand side, just by associativity of 
     composition. (This is essentially the same argument as the 
     argument given previously for interchangability of columns in 
     bi-rapports (= cross ratios), and it generalizes to 
     ``multi-rapports'', as considered in \cite{DL}.) 
     
     For 
     tri-rapports, it is not 
     true that the two {\em rows} of the matrix can be interchanged.

     The following equation is trivial, by repeated use of 
     (\ref{one}):
     \begin{equation}\left[ \begin{array}{ccc}A&B&C\\ D&E&F 
     \end{array}\right] ^{-1} =
	     \left[ \begin{array}{ccc}A&C&B\\ F&E&D 
	 \end{array}\right]\label{invv}\end{equation}
     
     \medskip
     Not all tri-rapports can be expressed as bi-rapports with the 
     same entries, but every 
     bi-rapport can be expressed as a tri-rapport: 
     \begin{prop}We have
     \begin{equation}\left[\begin{array}{cc}A&B\\
     C&D \end{array}\right] = \left[ \begin{array}{ccc}A&C&D\\
     B&A&B\end{array} 
     \right].\label{bi-tri}\end{equation}\label{BI-TRI}\end{prop}
     {\bf Proof.} This is just a re-interpretation of the commutative 
     diagram (\ref{hex1}); 
     Here, the triangle commutes by (\ref{one}). Hence the inner 
     pentagon commutes. The upper composite in it 
     is the bi-rapport considered; the lower composite is the 
     tri-rapport considered.
     
          \medskip
     
	We rewrite the classical ``cross ratio'' list, augmenting it with the 
	expression of the respective cross ratios (= bi-rapports) in terms of 
	tri-rapports, using Proposition \ref{BI-TRI}:
	
	\begin{align}\mu &=\left[ \begin{array}{cc} A&B\\ C&D
    \end{array}\right]&= \left[\begin{array}{ccc}A&C&D\\ B&A&B 
\end{array}\right]\\
   \mu ^{-1}&=\left[ \begin{array}{cc}A&B\\D&C
    \end{array}\right]&= \left[\begin{array}{ccc}A&D&C\\ B&A&B 
\end{array}\right]\\
    1-\mu &=\left[ \begin{array}{cc}A&C\\B&D
    \end{array}\right]&= \left[\begin{array}{ccc}A&B&D\\ C&A&C 
\end{array}\right]\label{PHI2}\\
    (1-\mu )^{-1}&=\left[ \begin{array}{cc}A&C\\D&B 
    \end{array}\right]&= \left[\begin{array}{ccc}A&D&B\\ C&A&C 
\end{array}\right]\\
    1-\mu ^{-1}&=\left[ \begin{array}{cc}A&D\\B&C
    \end{array}\right] &= \left[\begin{array}{ccc}A&B&C\\ D&A&D 
\end{array}\right]\\
    (1-\mu ^{-1})^{-1}&=\left[ \begin{array}{cc}A&D\\C&B
    \end{array}\right]&= \left[\begin{array}{ccc}A&C&B\\ D&A&D 
\end{array}\right]
    \end{align}

     \medskip
     
          With the $-1$ available as a ``uniform'' scalar, the six 
     $\mu$-expressions in the ``classical list'' above may be augmented by the six further 
     ones, obtained by 
     putting  minus sign on the right hand sides. The 
      scalars thus defined cannot in general be expressed as 
     cross ratios (bi-rapport) of four points, but can, by Proposition 
     \ref{BI-TRI} be expressed 
     as {\em tri-rapports} of four points. First, we have

     \medskip

     \begin{prop}We have
	 $$-\left[\begin{array}{cc}A&B\\ C&D \end{array}\right] = \left[ 
	 \begin{array}{ccc}A&B&D\\ C&A&B \end{array}\right] =
	     \left[ 
	 \begin{array}{ccc}A&C&B\\ B&A&D \end{array}\right]$$
	 \label{MINUS}\end{prop}
     {\bf Proof.} To prove the first equality,
      consider the diagram
     $$
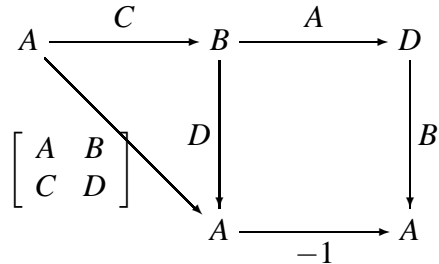
\begin{diagram}[nohug]A&\rTo^{C}&B&\rTo^{A}&D\\
     &\rdTo _{\left[\begin{array}{cc}A&B\\ C&D \end{array}\right] }&\dTo^{D}&&\dTo_{B}\\
     &&A&\rTo_{-1}&A
     \end{diagram}$$
     The triangle commutes, by definition of $(A,B;C,D)$; the 
     square commutes by (a variation of) the definition of $-1$. The clockwise composite 
     is the tri-rapport $(A,B,D;C,A,B)$. So we get that this 
     tri-rapport is $(A,B;C,D) \cdot (-1)$, proving the first 
     equality.
     Next, consider the diagram
     $$
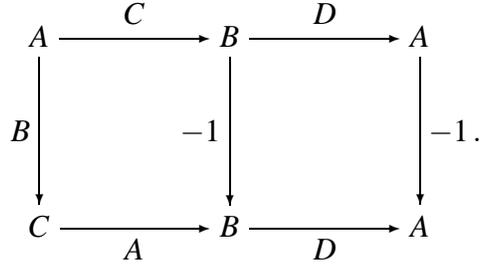
\begin{diagram}[nohug] A&\rTo ^{C}&B&\rTo^{D}&A\\
     \dTo^{B}&&\dTo^{-1}&&\dTo _{-1}\\
     C&\rTo_{A}&B&\rTo_{D}&A
     \end{diagram}.$$
     \begin{sloppypar}\noindent The clockwise composite is again 
     $-(A,B;C,D)$, the counterclockwise  is 
     $(A,C,B;B,A,D)$.\end{sloppypar}

\medskip

    Having the expressions in right hand column of the above table, we can give 
    tri-rapport expressions for the ``additive inverses'' of the six 
    scalars listed, using Proposition \ref{MINUS} and substitution 
    instances thereof. We refrain from using 
    arithmetic reductions like $-(1-\mu )=\mu -1$, because validity 
    of arithmetic 
    has not been assumed. We do, however, implicitly use that the 
    involutions $x\mapsto x^{-1}$ and $x\mapsto -x$ ($:= (-1)\cdot x)$ 
    do commute; this follows from $(-1)^{-1 }=-1$ (Proposition \ref{xx}).
    
    	\begin{align}-\mu &= \left[\begin{array}{ccc}A&B&D\\ C&A&B \end{array}\right]
	    &=\left[\begin{array}{ccc}A&C&B\\B&A&D \end{array}\right]\\
   -\mu ^{-1}&=\left[\begin{array}{ccc}A&B&C\\ D&A&B \end{array}\right]&=\left[\begin{array}{ccc}A&D&B\\B&A&C \end{array}\right]\\
    -(1-\mu ) &= \left[\begin{array}{ccc}A&C&D\\ B&A&C \end{array}\right]
    &=\left[\begin{array}{ccc}A&B&C\\C&A&D 
\end{array}\right]\label{min3}\\
    -(1-\mu )^{-1}&= \left[\begin{array}{ccc}A&C&B\\ D&A&C \end{array}\right]&=\left[\begin{array}{ccc}A&D&C\\C&A&B \end{array}\right]\\
    -(1-\mu ^{-1}) &= \left[\begin{array}{ccc}A&D&C\\ B&A&D \end{array}\right]&=\left[\begin{array}{ccc}A&B&D\\  D&A&C \end{array}\right]\\
    -(1-\mu ^{-1})^{-1}&= \left[\begin{array}{ccc}A&D&B\\ B&A&D \end{array}\right]&=\left[\begin{array}{ccc}A&C&D\\D&A&B \end{array}\right]
    \end{align}

    \medskip	

	\noindent{\bf Remark.} The classical way of dealing with the scalar $-1$, here 
	defined as a tri-rapport, is in terms of {\em harmonic 
	conjugates}:	         
  Given $A,B,C, H$. Then 
$$(A,B;C,H)= -1 $$ iff $H:B\to A$ equals the composite
$$\begin{diagram}
B&\rTo^{A}&C&\rTo^{B}&A\end{diagram}.$$
For, precomposing the composite with $C:A\to B$ gives $-1$, and 
 precomposing $H:B\to A$ with $C:A\to B$ gives $-1$ iff  $(A,B;C,H)=-1$.
 The classical  way of formulating this 
characterizing property of $H$  is: $H$ is the 
{\em harmonic conjugate} of $C$ w.r.to $A,B$.
    
    \medskip
The cross-ratios (bi-rapport) and the particular kind of tri-rapport 
considered in (\ref{hex2}) together equip each $K=\h (A,A) \cup \{0\}$ 
with enough structure for a {\em field} (provided 
sufficient equations can be secured), namely
\begin{itemize}
    \item the groupoid structure assumed for ${\LL }$ 
	gives the 
	multiplication (together with $0\cdot x = 0$ for all $x$).
	\item the cross ratio relation $(A,B;D,C)= (A,B;C,D)^{-1}$ 
	gives the multiplicative inversion (which anyway was given 
	apriori, since every arrow in a groupoid does have an inverse).
	\item the involution
	$(A,B;C,D) \mapsto (A,C;B,D)$ gives the (candidate for) 
	$x\mapsto 1-x$
	\item the tri-rapport considered in (\ref{hex2}) gives the 
	(candidate for) $-1$.
	\end{itemize}
	Then the addition $+$ may be defined by
	$$x+y:= x\cdot (1-\bigl( (-1)\cdot x^{-1}\cdot y \bigr),$$
	(together with $0+x=x$).
	
	\medskip
	
	We can now state the Theorem. We are 
	assuming a projective line candidate $\LL$, with commutative vertex 
	groups $\h (A,A)$, satisfying  (\ref{one}), 
	(\ref{two}) and the two hexagon conditions (\ref{hex1}) and 
	(\ref{hex2}) (these conditions are purely equational), 
	as well as the condition (\ref{AS}).
	 Assume  that $\LL$ satisfies these 
	 conditions, and assume finally 
	    that $K=\h (A,A) \cup \{0\}$ carries a field structure, 
	    with the field multiplication in $K^{*}$ (the group of multiplicative 
	    units of $K$) equal to the groupoid 
	    composition in $\h (A,A)$, and 
	    such that the operation $x \mapsto (1-x )$ (as 
	    given by (\ref{PHI2})) equals the operation $x\mapsto 1-x$ as given 
	    by the field structure. (Such a structure is unique, if it 
	    exists, since we argued that the addition $+$ is 
	    determined by the remaining projective-line-candidate 
	    structure. So the final assumption may possibly be satisfied 
	    automatically: it is a matter of the associative law for 
	    $+$, and of the distributive law.)
	     \begin{thm} Under these circumstances,  $\LL$ 
	    is isomorphic to the projective line $P(K^{2})$. More 
	    precisely, given three distinct points in $\LL$, then 
	    there is a unique isomorphism of projective line 
	    (-candidates) taking the three given points to $[0:1]$, $[1:0]$ 
	    and $[1:1]$, respectively.
	    \end{thm}
	    {\bf Proof.} Each $\h (A,A)$ is by construction of 
	    $K$ identified with $K^{*}$ ; and this identification is compatible with 
	    $x\mapsto 1-x$, by assumption. Then the result is a Corollary of the 
	    ``Fundamental Theorem'' for abstract projective lines 
	    over $K$, as formulated in \cite{APL} \S 3.

	    \medskip

    \noindent{\bf Appendix}

\medskip

\noindent The ``tri-rapport'' table has geometric content in the sense that it 
gives recipes for {\em geometric} construction of certain algebraic 
combinations. As an illustration, we give a geometric (tri-rappport) construction of the scalar $-3$ ($ = 
-(1-(-2))$) in terms of the scalar $-2$, (presented in terms of a 
bi-rapport).

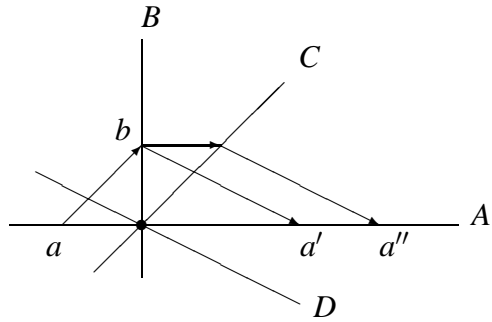
\begin{picture}(50,130)(-110,-40)

	\put(0,-20){\line(0,1){90}}
	\put(-18,-18){\line(1,1){72}}
	\put(-50,0){\line(1,0){170}}
	\put(-40,20){\line(2,-1){100}}
	\put(-30,0){\vector(1,1){30}}
	\put(0,30){\vector(2,-1){60}}
	\put(0,30){\vector(1,0){30}}
	\put(30,30){\vector(2,-1){60}}
	\put(125,0){$A$}
	\put(0,75){$B$}
	\put(60,60){$C$}
	\put(65,-35){$D$}
	\put(-36,-12){$a$}
	\put(60,-12){$a'$}
	\put(90,-12){$a''$}
	\put(-10, 32){$b$}
	\put(0,0){\circle*{4}}
	\end{picture}

\noindent We take the groupoid $\LL$ to have as objects the lines through a given 
point in a plane (indicated by a dot in the figure) (this is a 
standard representation of a projective line: as (unoriented) {\em directions} in 
a plane);   the arrows $A\to B$  
are the bijective linear maps between these lines 
(viewing them as 1-dimensional vector spaces, with the dot as zero). 
If $A\neq B$ (as in the figure), these linear maps are given by 
projection in a specific direction; thus $C:A\to B$ maps the point $a 
\in A$ to the point $b\in B$.
The cross ratio $(A,B;C,D)$  is a linear endo-map of the line 
$A$, and it takes $a$ to $a'$; 
this linear endo-map looks like it is ``multiplication by the 
scalar $-2$''. The tri-rapport $(A,B,C;C,A,D)$ takes $a$ to $a''$, and 
this looks like it is ``multiplication by the scalar $-3$'', in 
agreement with the tri-rapport formula
for $-(1-\mu)$ ($= \mu -1$) in 
(\ref{min3}).

\
    
    \noindent March 2010
    
\end{document}